\documentclass{amsart}
\usepackage{amssymb}
\newtheorem{thm}{Theorem}[section]
\newtheorem{lem}[thm]{Lemma}
\newtheorem{dn}[thm]{Definition}

\newtheorem{pro}[thm]{Proposition}

\def\rref#1{(\ref{#1})}

\newcommand{\daysub}[2][n]{\mbox{$#2_1,#2_2,\ldots ,#2_{#1}$}}
\def\bba{\begin{eqnarray}}
\def\eea{\end{eqnarray}}
\def\eee{\rule{.75ex}{1.5ex}\vskip1ex}
\def\eqeee{\mbox{ \rule{.75ex}{1.5ex}}}
\newcommand{\BBb}[1]{{\Bbb #1}}
\def\bK{{\Bbb K}}
\def\bC{{\Bbb C}}
\newcommand{\bbas}{\begin{eqnarray*}}
\newcommand{\eeas}{\end{eqnarray*}}
\newcommand{\bbar}{\begin{array}}
\newcommand{\eear}{\end{array}}
\newcommand{\bbs}{\begin{displaymath}}
\newcommand{\ees}{\end{displaymath}}
\newcommand{\bb}{\begin{equation}}
\newcommand{\eqbb}{\begin{equation}}
\newcommand{\va}{\varepsilon}
\def\lora{\longrightarrow}
\def\End{\mbox{\rm End}}
\newcommand{\tr}{\mbox{\rm tr}}
\def\loma{\longmapsto}
\def\ot{\otimes}
 \newcommand{\chH}{\mbox{$\check{H}$}}
 \newcommand{\chV}{\mbox{$\check{V}$}}
 \newcommand{\Hsq}{\mbox{$H^{\square} $}}
 \newcommand{\HLh}{\mbox{$H_{L^2}$}}
\newcommand{\Hcs}{\mbox{$H_{C^*}$}}
\newcommand{\B}{{\mathcal B}}
\renewcommand{\H}{{\mathcal H}}
\newsymbol\injarrow 131A
\newsymbol\surarrow 1310
\newcommand{\K}{{\Bbb K}}
\def\F{\mbox{$\mathcal F$}}
\def\G{\mbox{$\mathcal G$}}

\let\testint\int
\def\int{\displaystyle\testint}
\let\testiint\iint
\def\iint{\displaystyle\testiint}
\let\testiiint\iiint
\def\iiint{\displaystyle\testiiint}

\def\id{{\mathchoice{\mbox{\rm id}}
                    {\mbox{\rm id}}
                    {\mbox{\scriptsize\rm id}}
                    {\mbox{\tiny\rm id}} }}

\def\trace{\mbox{trace}}
\def\lam{\lambda}

\def\x{{\bf x}}
\def\u{{\bf u}}
\def\U{{\bf U}}
\def\bfxi{{\mathbf \xi}}
\def\bfeta{{\mathbf \eta}}
\def\w{{\bf w}}
\def\W{{\bf W}}

\def\mytitle{COMPACT QUANTUM GROUP $C^*$-ALGEBRAS AS\\ HOPF ALGEBRAS WITH APPROXIMATE UNIT }

\def\myperaddress{Institute of Mathematics, National Center for Science and Technology, P.O.Box 631, Bo Ho, VN-10.000, Hanoi, Vietnam}
\def\mythanks{DND and PHH would like to thank the Abdus Salam International Centre for Theoretical Physics for its hospitality and financial support.}

\begin{document}
\def\amshead{\date{April 30, 1999}
\title[Compact quantum groups $C^*$-algebras]{\mytitle}
\author{Do Ngoc Diep}
\author{Phung Ho Hai}
\author{Aderemi O. Kuku}
\address[Phung Ho Hai, Do Ngoc Diep]{\myperaddress}
\email[Do Ngoc Diep]{dndiep@ioit.ncst.ac.vn}
\email[Phung Ho Hai]{phung@ioit.ncst.ac.vn}
 \address[Aderemi O. Kuku]{International Center for Theoretical Physics\\ P.O. Box 586, 34100 Trieste, Italy.}
\email{kuku@ictp.trieste.it}
 \thanks{\mythanks}
 \subjclass{}
 \maketitle }

 \bibliographystyle{plain}

 \amshead

\begin{abstract} In this paper, we construct and study the representation theory of a Hopf $C^*$-algebra with approximate unit, which constitutes quantum analogue of a compact group $C^*$-algebra.
The construction is done by first introducing a convolution-product on an arbitrary Hopf algebra $H$ with integral, and then constructing the $L_2$ and $C^*$-envelopes of $H$ (with the new convolution-product) when $H$ is a compact Hopf $*$-algebra.
\end{abstract}

\section{introduction}
Compact quantum groups were introduced by Woronowicz \cite{woron1} and studied by several other authors \cite{andr1,ed-ru1,daele1}, as non-commutative analogues of the function algebras on compact groups. The aim of this work is to construct some quantum analogues of the compact group $C^*$-algebras, which turn out to be non-cocommutative Hopf $C^*$-algebras with approximate unit.

More precisely, starting from a compact Hopf $*$-algebra, we construct a Hopf $C^*$-algebra with approximate unit in such a way that when the original compact Hopf $*$-algebra is the function algebra on certain compact group $G$, the resulting Hopf $C^*$-algebra reduces to the group $C^*$-algebra $C^*(G)$.

Recall that for a compact group $G$, the group $C^*$-algebra $C^*(G)$ is the completion by operator norm of the algebra $L^1(G)$ of absolutely integrable complex-valued functions on $G$. The product on $L^1(G)$ is the convolution-product, given by
$$(g*f)(x)=\int_{G}f(y)g(y^{-1}x)dy.$$
With respect to this product, $C^*(G)$ is a non-unital $C^*$-algebra. In that case, there is an important notion of $\delta$-type sequences in $C^*(G)$, that approximate the unity. Such a sequence is called an approximate unit in $C^*(G)$. Moreover, the action of the group $G$ on any representation can be recovered from the corresponding action of the algebra $C^*(G)$, by using $\delta$-type sequences. In our case, however, we do not know what the quantum group defining our $C^*$-algebra is, hence our defining a Hopf algebra structure to exhibit the ``group property'' of our algebra.

Our construction is done in several steps. First we construct a convolution product on an arbitrary Hopf algebras with integrals, being motivated by the classical convolution product for $L_1(G)$, the completion of which (under the operator norm) defines $C^*(G)$ (see \ref{chH}). Denoting by \chH\ the vector space $H$ with the new convolution product, we show that \chH\ is a non-unital algebra that is an ideal of $H^*$ (the dual of $H$); that the category of completely reducible $H$-comodules is equivalent to the category of completely reducible \chH-comodules (\ref{thmact-coact}, \ref{chh-hsq}) and that when $H$ is co-semi-simple then \chH\ is isomorphic to a direct sum of full endomorphism rings of simple $H$-comodules \rref{thmchH}.

In section \ref{sec3}, we now focus on the case when $H$ is a compact Hopf $*$-algebra. Through the works of Woronowicz, Koorwinder and Djikhuizen, \cite{woron1,dk}, compact Hopf $*$-algebra are known to play a role analogous to that of the algebra of functions on compact groups in the classical theory. Since by definition, $H$ has a complex scalar product, we complete it to a Hilbert space \HLh, which is the quantum analogue of the algebra $L^2(G)$. We show in \ref{extend-pr} that the convolution product on \chH\ extends to \HLh\ and that \HLh\ possesses a (topological) coproduct and an approximate antipode and hence is a Hopf algebra with approximate unit (\ref{thm-hlh}). Since the antipode on $H$ is not involutive, the $*$-structre defined on $H$ by $h^\star:=S(h^*)$ cannot be extended to \HLh\ because it is not continous with respect to the given norm. To overcome this problem, we pass to the $C^*$-enveloping of \HLh\ to obtain \Hcs\ that is the required Hopf $C^*$-algebra with approximate unit (see 4.5). We show that there is a one-one correspondence between the irreducible unitary representations of \Hcs\ and those of \chH\ (see 4.4).

 \section{The Convolution Product and Co-semi-simplicity}\label{sec1}
We work over an algebraically closed field $\bK$ of characteristic zero.
 Let $(H,m,\eta)$ be an algebra over a field $\K$, where $m$ denotes the product, $\eta$ denotes the map $\K\lora H$, $1_{\K}\lora 1_H$. A bialgebra structure on $H$ is a pair of linear maps $\Delta:H\lora H\ot H,\ \va:H\lora \K$, satisfying 
 \begin{itemize}\item $\Delta$ and $\va$ are homomorphism of algebras.
 \item $(\va\ot \id_H)\Delta= (\id_H\ot\va)\Delta=\id_H.$
 \item $(\Delta\ot\id_H)\Delta=(\id_H\ot \Delta)\Delta.$
 \end{itemize}
 An antipode on $H$ is a linear map $H\lora H$, satisfying
 \begin{itemize}
 \item $m\circ(S\ot\id_H)\Delta=m\circ(\id_H\ot S)\Delta=\eta\va.$
 \end{itemize}
 A bialgebra equipped with an antipode is called Hopf algebra. The antipode is then uniquely determined.

 Let $(H,m,\eta,\Delta,\va)$ be a Hopf algebra. A right coaction of $H$ on a vector space $V$ is a linear map $\delta:V\lora V\ot H$, satisfying
 \begin{equation}\label{coaction}
 \bbar{c}(\id_V\ot\va)\delta=\id_V,\\
 (\delta\ot\id_H)\delta=(\id_V\ot \Delta)\delta,\eear\end{equation}
in the first identity we identify $V$ with $V \otimes \mathbb{C}$. $V$ is then called a right $H$-co-module.
 
We shall frequently use Sweedler's notation, in particular $\Delta(x):=\sum_{(x)}x_{(1)}\ot x_{(2)}=x_1\ot x_2, \delta(v)=\sum_{(v)}v_{(0)}\ot v_{(1)}=v_0\ot v_1.$

 The elements $v_1$'s in the presentation $\delta(v)=v_0\ot v_1$ are called coefficients of the coaction $\delta$, and the space they span is called coefficient space. This space is a subcoalgebra of $H$. To see this, fix a basis $\daysub[d]{x}$. Then, the coaction is given by $\delta(x_i)=x_j\ot a_i^j$ and the $\{a_i^j\}$ span the coefficient space. On the other hand, from \rref{coaction} it follows $\Delta(a_i^j)=a^j_k\ot a^k_i$, $\va(a_i^k)=\delta^k_i$.
 The comodule $V$ is simple iff $\{a_i^j\}$ is a basis for the coefficient space, see \cite{schmuedgen1}.

 A left (right) integral on a Hopf algebra is a (non-trivial) linear functional $\int:H\lora \K$, which is a left (right) $H$-comodule homomorphism, where $H$ is a left (right) $H$-comodule by means of the coproduct $\Delta$ and $\K$ is a left (right) $H$-comodule by means of the unity map $\eta$. Explicitly, a left (resp. right) integral $\int_l$ $\left(\mbox{resp.,} \int_r\right)$ on $H$ satisfies
\bba\label{axiom-int}\int_l(x)=x_1\int_l(x_2)\quad \left(\mbox{resp.} \int_r(x)=\int_r(x_1)\cdot x_2\right).\eea

It was shown by Sullivan \cite{sul} that the integral on a Hopf algebra, if it exists, is defined uniquely up to a constant.  Further, we have
\begin{lem}\label{lemnon-degenerate}\cite{sweedler,lin1,stefan1} Let $\int $ be a left integral on $H$. Then the bilinear form $b(g,h)=\int(gS(h))$ is non-degenerate on $H$, that is
\bba\label{non-degenerate}\bbar{l} \int(gS(h))=0, \forall h\quad \Longrightarrow\quad g=0.\\
\int(gS(h))=0, \forall g\quad \Longrightarrow\quad h=0.\eear\eea\end{lem}

In the rest of this work, we fix a left integral $\int$ on $H$.

We mention an identity, due originally to Sweedler \cite{sul}, which plays a crucial role in our computations:
\begin{equation}\label{sweedler-id}
 \int(gS(h_1))\cdot  h_2=g_1\int(g_2S(h)),\end{equation}
if $\Delta(g) = g_1 \otimes g_2$ and $\Delta(h) = h_1 \otimes h_2$ in the notations of Sweedler.

The coalgebra structure on $H$ induces an algebra structure on its dual $H^*$ -- the space of linear forms on $H$. The product of $\phi,\psi$ in $H^*$ is given by
$$\phi*\psi(h)=\phi(h_1)\psi(h_2), \forall h\in H.$$
The unit element in $H^*$ is the counit $\va$ of $H$.

 \subsection*{Convolution product on $H$}

 We first recall the classical structure. Let $G$ be a compact group. Then there exists a unique normalized Haar measure on $G$ which induces the Haar integral on $L^1(G)$. The $*$-product is defined on $L^1(G)$ as follows:
 \begin{equation}\label{conv-pr}
 (g*f)(x):=\int_{G}f(y)g(y^{-1}x)dy.\end{equation}

 Now, let $H$ be a Hopf algebra with an integral. Being motivated by \rref{conv-pr} we define the convolution product on $H$ by:
 \begin{equation}\label{conv-pr1}
 g*f:=\int (fS(g_1))\cdot  g_2.\end{equation}
According to \rref{sweedler-id}, we also have $g*f=f_1\int(f_2S(g)).$
 \begin{lem}\label{chH}
 $H$ equipped with the $*$-product defined above is a (non-unital) algebra.\end{lem}
 \proof
 We only have to check the associativity: 
 \begin{eqnarray*}
 (h*g)*f&=& \int\left(fS((h*g)_1)\right)\cdot  (h*g)_2\\
  &=& \int\left(fS\left(\int(gS(h_1)\cdot  h_2)\right)\right)\cdot  h_2\\
 &=& \int\left(fS(h_2)\cdot \int(gS(h_1)\right)\cdot  h_3\\
 &=& \int\left(f\int(gS(h_1))\cdot  S(h_2)\right)\cdot  h_3 \mbox{(using \rref{sweedler-id})}\\
 &=&\int\left(fS(g_1)\int(g_2S(h_1))\right)\cdot  h_2\\
 &=&\int\left(fS(g_1)\cdot \int(g_2S(h_1))\right)\cdot  h_2\\
 &=& \int\left(\int(fS(g_1))\cdot  g_2 S(h_1)\right)\cdot  h_2\\
 &=& \int((g*f)S(h_1))\cdot  h_2\\
 &=& h*(g*f).\end{eqnarray*}
 Lemma \ref{chH} is proved. \eee

We denote by $\check H$ the vector space $H$, equipped with the convolution product $*$.

\subsection*{$H$-comodules and \chH-modules}
 We now study the correspondence between $H$-comodules and \chH-modules.
 Let $V$ be an $H$-comodule. We define an action of $\chH$ on $V$ as follows:
 \begin{equation}\label{coact-act}\textstyle
 h*v:=v_0\int( v_1S(h))\end{equation}
 We check the associativity: 
 \bbas g*(h*v)&=& v_0\int(v_1S(g))\cdot \int(v_2S(h))\\
(\mbox{using }\rref{sweedler-id}) &=& v_0\int\left(v_1S(h_1)\int(h_2S(g))\right)\\
&=& v_0\int(v_1S(g*h))\\
 &=& (g*h)*v.\eeas

 Let $\phi:V\lora W$ be a morphism of $H$-comodules, i.e., $\phi(v)_0\ot \phi(v)_1=\phi(v_0)\ot v_1$. Then, for $h\in H, v\in v$,
 \bbas \phi(h*v)&=& \phi(v_0\int(v_1S(h))\\
 &=& \phi(v_0)\int(v_1S(h))\\
 &=& \phi(v)_0\int(\phi(v)_1S(h))\\
 &=&h*\phi(v).\eeas

 Thus, $\phi$ is a morphism of \chH\ modules. We therefore have a functor $\mathcal F$ from the category of $H$-comodules into the category of \chH-modules, which is the identity functor on the underlying category of vector spaces.
\begin{pro}\label{funF} The functor $\mathcal F$ defined above is full, faithful and exact.\end{pro}
\proof From the definition of $\mathcal F$, we see that, as vector spaces, $\F(V)=V$ and $\F(\phi)=\phi$. Thus, the functor $\F$ is faithful and exact. It remains to show that $\F$ is full, which amounts to showing that if $\phi:V\lora W$ is a morphism of $\chH$-module then it is a morphism of $H$-comodule. By assumption, $\phi$  satisfies
$$ h*(\phi(v))=\phi(h*v),\quad \forall h\in\chH$$
or, equivalently,
$$\phi(v)_0\int\left(\phi(v_1)S(h)\right)=\phi(v_0)\int(\left(v_1S(h)\right), \quad \forall h\in\chH.$$
Since $\int$ is faithful (see \rref{non-degenerate}), we conclude that 
$$\phi(v_0)\ot \phi(v_1)=\phi(v_0)\ot v_1.$$
In other words, $\phi$ is a homomorphism of $H$-comodules.\eee

 Let $V$ now be a cyclic \chH-module, i.e., there exists an element $\bar{v}\in V$, such that all $v\in V$ is obtained from $\bar{v}$ by the action of some $f\in H$. We then define the coaction of $H$ on $V$ by 
 \begin{equation}\label{act-coact}
 \delta(v):=f_1* \bar{v}\ot S(f_2),\end{equation}
where $f$ is such that $f*\bar{v}=v$.

 First, we have to show that this coaction is well defined, which means that it does not depend on the choice of the representative elements $\bar{v}$ and $f$. To show that the definition does not depend on $\bar{v}$ we let $\bar{v}=g*\tilde{v}$ and show that the definition does not change when $\bar v$ is replaced by $\tilde{v}$, which means
\bbs f_1*\bar{v}\ot S(f_2)=(f*g)_1*\tilde{v}\ot S((f*g)_2).\ees
Replacing $\bar v$ on the left-hand side of this equation by $g*\tilde{v}$ and canceling $\tilde{v}$ in both sides, one is led to the following equation
\bbs f_1*g\ot f_2=(f*g)_1\ot (f*g)_2,\ees
which follows immediately from the definition.

To show the independence on the choice of $f$ we assume $f*\bar{v}=0$ and show that $f_1*\bar{v}\ot S(f_2)=0$. Indeed, we have 
$$\int(f_2S(g))\cdot  f_1*\bar{v}=(g*f)*\bar{v}=g*(f*v)=0,$$ for all $g$. Hence, according to Lemma \ref{lemnon-degenerate}  $f_1*\bar{v}\ot f_2=0$.

 We now proceed to check the co-associativity and co-unitary:
 \bbas (\id_H\ot\delta)\delta(v)&=&\delta(f_1*\bar{v})\ot S(f_2)\\
 &=& f_1*\bar{v}\ot f_2\ot f_3)\\
 &=& f_1*\bar{v}\ot \Delta(f_2)\\
 &=&(\id\ot\Delta_H)\delta(v).\eeas
 \bbas(\id_V\ot \va)\delta f(v)&=&f_1*\bar{v}\ot\va(f_2)\\
 &=&f*\bar{v}=v.
 \eeas

Since simple modules are cyclic, we obtain a functor $\mathcal G$ from the category of completely reducible \chH-modules to the category of $H$-comodules.

 \begin{thm}\label{thmact-coact}
 The category of completely reducible $H$-comodules is equivalent to the category of completely reducible \chH-modules. 
\end{thm}
\proof We show that the functors $\mathcal F\circ \mathcal G$ and $\mathcal G\circ \mathcal F$ are identity-functors on the categories of completely reducible $H$-comodules and \chH-comodules, respectively. 

Let $V$ be a simple $H$-comodule. Then $\chV={\mathcal F}(V)$ is a simple \chH-module. Let us fix $\bar{v}\in\chV$ and for $v\in\chV$, let $f\in\chH$ be such that $f*\bar{v}=v$. By definition, we have
\bbs v=f*\bar{v}=\bar{v}_0\int(\bar{v}_1S(f)).\ees
The coaction of $H$ on $\G\circ \F(V)$ is 
\bbas\delta_{\mathcal{G}\circ\mathcal{F}(V)}(v)&=&f_1*\bar{v}\ot f_2\\
&=&\bar{v}_0\int(\bar{v}_1S(f_1))\ot f_2\\
&=&\bar v_0\ot \bar v_1\int(v_2S(f))\\
&=&\delta_V(v).\eeas
Thus, $\G\circ\F$ is the identity functor. The assertion for $\F\circ\G$ is proved analogously. \eee

 \subsection*{\chH\ is an ideal of $H^*$}
 In the previous section we have seen that there exists a correspondence between $H$-comodules and \chH-modules. On the other hand, there exits a one-to-one correspondence between $H$-comodules and rational $H^*$-modules. It is then natural to ask about the relationship between $H^*$ and \chH. We now show that \chH\ is isomorphic to the left ideal generated by the integrals in $H^*$. 

 It is well-known, that the rational  submodule  \Hsq\ of $H^*$, considered as a left module on itself, i.e. the sum of all left ideals of $H^*$, which are finite dimensional (over $\bK$), is an $H$-Hopf module (with an appropriate $H$-action) and hence isomorphic to the tensor product of $H$ with the space spanned by the integrals \cite[Thm 5.1.3]{sweedler}. Since the space of integrals is one-dimensional, we have an isomorphism between the two vector spaces $H$ and \Hsq. This isomorphism can be given explicitly as follows.

 By means of the integral, every element of $H$ can by considered as a linear functional on $H$ itself: $H\ni h\loma \int^h\in H^*: \int^h(g):=\int(gS(h)).$ 

 \begin{pro}\label{chh-hsq} The map $H\ni h\loma \int^h\in H^*$ defined above is an isomorphism of algebras $\chH\lora H^\square\subset H^*$.\end{pro}
 \proof
 We have
 \bbas \int^{f*h}(g)&=& \int\left(gS(h_1\int (h_2Sf))\right)\\
 &=&\int(gS(h_1))\int(h_2S(f))\\
 &=&\int(g_1S(f))\int(g_2S(h))\\
 &=&(\int^f*\int^h)(g).\nonumber \eqeee\eeas
 
 In the strictly algebraic sense, \chH\ is not a Hopf algebra, unless $H$ is finite-dimensional. In the next section we will show that for compact Hopf $*$-algebra, there exists a natural topology in $H$ such that the completion of $H$ with respect to this topology is a topological Hopf algebra with approximate unit.

\subsection*{Co-semi-simple Hopf algebras}
A Hopf algebra $H$ is called co-semisimple if any finite dimensional $H$-comodule decomposes into a direct sum of simple comodules. A Hopf algebra is co-semisimple if and only if it possesses an integral whose value at its unit element is nonzero. In this case, left and right integrals are equal \cite{sweedler}.

 Let $H$ be co-semi-simple Hopf algebra. Then it decomposes into a direct sum of simple sub-coalgebras, each of which is the coefficient space of a simple $H$-comodule \cite{sweedler}
 \bba\label{coalg-decom}H\cong\bigoplus_{\lambda\in\Lambda}H_\lambda.\eea
  The set $\Lambda$ contains an element $0$ for which $H_0\cong\K$. The integral computed on $H_\lambda$ is zero for $\lambda\neq 0$.

We consider the discrete topology on $\Lambda$. Let $C_0(\Lambda)$ denote the set of all compact subsets in $\Lambda$ containing 0. For any compact $K$, let
 \bba\label{hk} H_K:=\bigoplus_{\lambda\in K}H_\lambda.\eea
 Then $H_K$ are subcoalgebra of $H$. For any $f\in H$, $f=\sum_{\lambda\in K} f_\lambda, f_\lambda\in H_\lambda,$ for some compact $K$.

 For each $\lambda\in\Lambda$, let $V_\lambda$ be the corresponding simple $H$-comodule. Note that the $V_\lambda$ are finite dimensional for all $\lambda\in \Lambda$.

 The isomorphism \rref{coalg-decom} becomes now an isomorphism of algebras between \chH\ and the direct sum of endomorphism ring of $V_\lambda$.
 \bba\label{alg-decom}
 \chH\cong\bigoplus_{\lambda\in\Lambda}\chH_\lambda\cong\bigoplus_{\lambda\in\Lambda}\End_\K(V_\lambda).\eea

Thus we have proved
 \begin{thm}\label{thmchH}
  Let $H$ be a co-semisimple Hopf algebra. Then \chH\ is isomorphic to the direct sum of full endomorphism rings of simple $H$-comodules.
 \end{thm}

 The algebra \chH\ does not have a unit element. Adding a unit to this algebra is problematic when we are dealing with the norm -- the unit element is some thing like `` the Dirac delta function " which never has finite norm. Instead we have a notion of $\delta$-type sequences, which approximate the unit.

 \begin{dn}\label{dnau}\rm 
\begin{itemize}\item[1)]
 Let $A$ be an algebra without unit. A system $\{e_i,i\in I\}$ of idempotents in $A$ is an {\it approximate unit} if
 \begin{itemize}\item[(i)] $I$ is a partially ordered set,
 \item[(ii)] For any $a\in A$, there exists $i=i(a)$, such that $e_ia=ae_j=a$, for all $j\geq i$.\end{itemize}
 For an approximate unit in a bialgebra we require further that
\begin{itemize} \item[(iii)] $\va(e_i)=1, \forall i\in I$.\end{itemize}
  For an algebra with involution we require that
\begin{itemize} \item[(iv)] There exist $f_i$, such that $e_i=f_if_i^*$, for all $i\in I$.
 \end{itemize}
\item[2)]For a topological algebra, the condition (ii) above is replaced by
\begin{itemize}\item[(ii')] the nets $\{ e_ia|i\in I\}$ and $\{ ae_i|i\in I\}$ converge to $a$.\end{itemize}
\item[3)]A Hopf algebra with approximate unit is a bialgebra with approximate unit together with a system of endomorphism $\{ S_i|i\in I\}$, called an {\it approximate antipode}, satisfying
$$m(S_i\ot \id)\Delta=m(\id\ot S_i)\Delta=e_i\va.$$
\end{itemize}
 \end{dn}

 The existence of such a sequence in our \chH\ is obvious. Indeed, let $e_K$ be the unit element in \chH$_K$. Then $\{e_K,K\in C_0(\Lambda)\}$ is an approximate unit in \chH. Thus \chH\ is an algebra with approximate unit.

 \section{Compact Hopf $*$-algebras}\label{sec2}

 In this section  we construct from a compact Hopf $*$-algebra a bialgebra which is a Hilbert space. Thus, $\Bbb K=\Bbb C$. This is the first step toward the construction of our Hopf $C^*$-algebra. A good reference on compact Hopf $*$-algebra, where the algebra is referred to as CQG-algebra, is Dijkhuizen and Koornwinder \cite{dk}.

  By definition, a compact Hopf $*$-algebra over $\Bbb C$ is a co-semisimple Hopf algebra with an involutive $\Bbb C$-anti-linear anti-homomorphism $^*$ such that every simple comodule is unitarizable, i.e. we can define a scalar product on this comodule such that
 \begin{equation}\label{unitary-mod} <v_0,w>S(v_1)=<v,w_0>w_1^*,\mbox{ for }v,w\in H.\end{equation}

  For any orthonormal basis of this comodule, the corresponding coefficient matrix satisfy the orthogonality condition. More precisely, let $\x_1,\x_2,\ldots, \x_d$ be an orthonormal basis of the comodule and $\U=(\u_i^j)$ be the corresponding coefficient matrix, i.e., $\delta(\x_i)=\x_k\ot \u_i^k$. Then $\U$ satisfies $\U \U^*=\U^*\U=I$, where ${\U^*}^j_i:=\u_j^{i*}$.

\begin{lem}\label{lemssquare}Let $H$ be a co-semisimple Hopf algebra. Then the square of the antipode on $H$ is co-inner, i.e., it can be given in terms of an invertible element $q$ of $H^*$ (i.e. a linear form of $H$):
$$S^2(h)=q(h_1)h_2q^{-1}(h_2),$$ where
the linear form $q$ is given by $q(h)=\int(S^2(h_1)S(h_2)).$
\end{lem}
\proof First we show that $S^2(h_1)q(h_2)=q(h_1)h_2,$ or equivalently 
\bba\label{sq}S^2(h_1)\int(S^2(h_2)S(h_3))=\int(S^2(h_1)S(h_2))\cdot h_3.\eea
We have
\bbas \int(S^2(h_2)S(h_3))\cdot h_4S(h_1)&=&S^2(h_2)\int(S^2(h_3)S(h_4))\cdot S(h_1) (\mbox{by } \rref{sweedler-id})\\
&=&\int(S^2(h_1)S(h_2)).\eeas

Thus, the lemma will be proved if we can show that $q$ is invertible as an element of $H^*$.

Let $V$ be a finite dimensional $H$-comodule and $V^{**}$ be its double-dual. As vector space, $V^{**}$ is isomorphic to $V$ and the coaction of $H$ on $V^{**}$ is given by $\delta_{V^{**}}(v)=v_0\ot S^2(v_1)$. Equation \rref{sq} shows that the map $v\loma v_0\ot q(v_1):V\lora V^{**}$ is a morphism of $H$-comodules. If $V$ is simple then $V^*$ and hence $V^{**}$ are simple. Therefore, the map above should be zero or invertible. To see that it cannot be zero we fix a basis $\x_1,\x_2,\ldots,\x_d$ of $V$ and let $\U=(\u_i^j)$ be the coefficient matrix. Since have $\Delta(\u_i^j)=\sum_k\u^j_k\ot \u_i^k, \va(\u_i^j)=\delta_i^j$,
$$q(\sum_i\u_i^i)=\int(\u^k_iS(\u_k^i))=d.$$

Thus, the map  $v\loma v_0\ot q(v_1):V\lora V^{**}$ is an  isomorphism of $H$-comodules. Therefore the form $q$ is invertible. \eee

Set $Q_i^j=q(\u_i^j)$. Then, according to Lemma \ref{lemssquare} the matrix $Q$ is the matrix of the isomorphism $V\lora V^{**}$ with respect to the basis $\x_1,\x_2,\ldots,\x_d$ of $V$ and we have
\begin{equation}\label{ssquare}S^2(\u_i^j)=Q^j_k\u^k_l{Q^{-1}}^l_j.\end{equation}

The matrix $Q$ is very important in the study of $H$ and will be called {\it reflection matrix}. If $V$ is irreducible then the integral on $\u_i^jS(\u_k^{l})$ can be given in terms of $Q$. In fact, from the left-invariance of $\int$ we have
 \begin{equation}\label{compute-int}
  \int(\u_i^jS(\u_k^l))=\u_m^jS(\u_k^n)\int(\u_i^mS(\u_n^l))
 \end{equation}
 Or equivalently
 \bbs \int(\u_i^jS(\u_k^l))\cdot  \u_n^k=\u_m^j\int(\u_i^mS(\u_n^l)).\ees
 Since $\{\u_i^j\}$ are linearly independent, $\int(\u_i^jS(\u_k^l))=\delta_k^jC_i^l$ for some matrix $C=(C_i^l)$. 

 On the other hand, according to \rref{ssquare}, we have $\u^m_iQ_m^nS(\u_n^j)=Q^j_i$. Substituting this into \rref{compute-int} we get $\delta_n^mC_i^jQ_m^n=Q^i_j$. Thus $\tr(Q)\neq 0$ and $C_i^j=Q^j_i/\tr(Q)$, and we get
 \begin{equation}\label{int}
 \int(\u_i^jS(\u^{l}_k))=\delta_k^jQ_i^l/\tr(Q).\end{equation}

 Analogously we have
 \begin{equation}\label{int2}
 \int(S(\u_i^{j})\u^l_k)=\delta_k^j(Q^{-1})_i^l/\tr(Q^{-1}).\end{equation}

  From Equations \rref{int} and \rref{int2} we get the rule for the convolution product for coefficients of simple comodules. Let $V_\lambda$ and $V_\mu$ be simple comodules and ${\U_\lambda},{\U_\mu}$ be the coefficient matrix with respect to some (orthogonal) bases of $V_\lambda$ and $V_\mu$, and $Q_\lambda:=q(\U_\lambda), Q_\mu:=q(\U_\mu)$. An element $f$ of $H_\lambda$ can be represented by a matrix $F_\lambda$: $f={F_\lambda}_j^i{\u_\lambda}^j_i=\trace(F_\lam \U_\lam).$ Let $c_\lam=\trace(C_\lam \U_\lam), d_\mu=\trace(D_\mu \U_\mu)$. The convolution product has the form
  \bba\label{product-rule} c_\lambda * d_\lambda&=&{C_\lam}^i_j{D_\lam}^k_l(\u_\lambda)_i^j*(\u_\mu)_k^l\nonumber\\
&=&{C_\lam}^i_j{D_\lam}^k_l\delta^\lambda_\mu \frac{{Q_\lambda}^j_k}{\tr(Q_\lambda)}(\u_\lambda)^l_i\nonumber\quad(\mbox{from \rref{conv-pr1} and \rref{int}})\\
&=&\frac{C_\lambda Q_\lambda D_\lambda}{\tr(Q)},\label{conv-pr2}\eea
 where in the left-hand side is the matrix product: $(CD)_i^j=C^i_kD^k_j$.
  Thus, with the $*$-product, $H_\lambda$ becomes a unital algebra, denoted by $\chH_\lambda$, and $\chH_\lam\cong\End_{\Bbb C}(V_\lambda).$ Moreover, for any compact $K\subset\Lambda$, $\chH_K:=\bigoplus_{\lambda\in K}\chH_\lambda$ is also a unital algebra.

 Let $V_\lambda$ be a simple $H$-comodule, $\chV_\lambda$ be the corresponding \chH-module. The representation of \chH\ on $\chV_\lambda$ will be denoted by $\pi_\lambda$.  Then according to \rref{coact-act}, the action of $f=\trace(F_\lam \U_\lam)$ on $\x_i$ is given by
 \begin{equation}\label{f-action}f*\x_i=\x_j\int(\u^j_iS(f))=\x_j{F_\lam}^j_l{Q_\lam}^l_i.\end{equation}
 In other words, the matrix of $\pi_\lambda(f)$ with respect to the basis $\{\x_i\}$ is $F_\lam Q_\lam$.

  For compact Hopf $*$-algebra, we can show that the matrix $Q$, with respect to an orthonormal basis, is positive definite \cite{dk,woron1}. Let $V^*$ be the dual to the comodule $V$. Let $\bfxi_1,\bfxi_2,\ldots,\bfxi_d$ be the basis of $V^*$, dual to the basis $\x_1,\x_2,\ldots,\x_d$. The corresponding coefficient matrix is then $S(\U^t)$.
 Since $V^*$ is unitarizable, there exists a basis $\bfeta_1,\bfeta_2,\ldots, \bfeta_d$ of $V^*$ which is orthonormal with respect to a scalar product, such that the corresponding coefficient matrix $\W=(\w_i^j)$ satisfies $S(\W)=\W^{*}$. Let $T$ be a matrix such that $\xi^i=\bfeta_jT_j^i$ then
  \begin{equation}\label{WtoU}\w^i_j={T^{-1}}^i_kS(\u_l^k)T^l_j={T^{-1}}^i_k\u^{l*}_kT^l_j,\end{equation}
  since $(\u^i_j)$ is also unitary. Therefore, by the involutivity of $*$
  \bbs S(\w_i^j)=\w^{i*}_j={\bar{T}}^{-1i}_{\:k}\u_k^l\bar{T}^l_j.\ees
  Substituting $S(w_i^j)$ into the preceding equation and using \rref{ssquare} we, get $Q=\mbox{const}\cdot  T^*T.$ A direct consequence of this fact is that $\int$ is positive definite on $H$. 

  For any $\lambda\in\Lambda$, let $\lambda^*$ be such that $(V_\lambda)^*=V_{\lambda^*}$. According to Lemma \ref{lemssquare}, $*$ is an involutive map on $\Lambda$. From above we see that the involution $*$ on $H$ maps $H_\lambda$ onto $H_{\lambda^*}$. Since the integral is zero on $H_\lambda$ except for $\lambda=0$, we see that
  $\int (a)=\overline{\int (a^*).}$
  Consequently 
 $\int$  defines a scalar product $<,>$ on $H$, $<a,b>:=\int(ab^*)$ giving rise to a norm on $H$ called $L^2$-norm. Let \HLh\ ddenote the completion of $H$ with respect to the $L^2$-norm. Then \HLh\ is a Hilbert space.

  \subsection*{The Hopf Algebra  \HLh}

  Our aim is to define a Hopf algebra structure on \HLh, where the product is an extension of the convolution product on \chH. First we have to extend the convolution product to \HLh. This can be done provided that we showed that this product is continuous with respect to the $L^2$-norm. 

  \begin{lem}\label{extend-pr}
  The convolution product on \chH\ ssatisfies
  $$\|f*g\|_{L^2}\leq \|f\|_{L^2}\cdot  \|g\|_{L^2}.$$\end{lem}
  \proof
  It is sufficient to show this inequality for $f$ and $g$
  belonging to the same coefficient space of a simple comodule. Then, in that 
  case, $f$ and $g$ can be represented in the form 
  $f = \tr(CU), \quad g = \tr(DU),$
  for a unitary multiplicative matrix $U$ and some matrices $C$, $D$ with
  complex scalar entries. Using \rref{int}, \rref{product-rule} the indicated inequality is
  equivalent to the following
  $$\frac{\tr(CQC^*)\tr(DQD^*)}{tr(Q)^2} \geq
  \frac{\tr(DQC Q (DQC)^*)}{\tr(Q)^3}.$$
 Let $T$ be such that $Q = TT^*$. Set $C_1 = CT,$ $D_1 = DT$
   the above inequality has the form
  $$\tr(C_1C_1^*)\tr(D_1D_1^*)\tr(TT^*) \geq \tr(C_1T^*D_1D_1^*TC_1^*).$$
  The last inequality follows immediately from the Minkowski inequality since
 $$\tr(CC^*)=\sum|c_i^j|^2.  \eqeee$$

  Therefore, the convolution product on \chH\ can be extended to \HLh.
  \begin{lem}
  The family $\{ e_K \vert K\in C_0(\Lambda) \}$, where $e_K$ is the unit element of $\chH_K$, is an approximate
  unit in $H_{L^2}$.
  \end{lem}
  \proof
  What we need to show is that for any $h\in H_{L^2}$, there exists a
  composition sequence $K_1 \subset K_2 \subset \dots $ such that
  $$\lim_{n\to \infty} \Vert e_{K_n}*h - h\Vert = \lim_{n\to \infty}
  \Vert h*e_{K_n} - h\Vert = 0$$
  Now, for any $h\in H_{L^2}$, there exists at most a countable set of $\lambda$,
  $\lambda \in \Lambda$ such that $h_\lambda = h*e_\lambda \ne 0$ and $h =
  \sum_\lambda h_\lambda$. The last series converges absolutely whence the
  assertion follows. \eee

   Then we have to define the coproduct. Notice that, the coproduct on \HLh, if it exists, should be dual to the original product on $H$ by means of the integral. Thus, we consider, for an element $h\in H$, a linear functional $\phi_h:H\ot H\lora \Bbb C$, $\phi_h(g\ot f)=\int(hgf)$. This is obviously continuous hence is extendable on $\HLh\hat{\ot}\HLh$ ($\hat{\ot}$ denotes the tensor product of Hilbert spaces). Hence, by Riesz theorem, there exists an element $\Delta_*(f)$ of $\HLh\hat{\ot}\HLh$, such that $\phi_h(g\ot f)=\int\int(\Delta_*(h),g\ot f)$. Thus, we get a map $H\lora \HLh$, $h\loma\Delta_*(h)$. Again this is continuous and hence induces a map $\HLh\lora \HLh\hat{\ot}\HLh$ which is the coproduct on \HLh.

  The counit is given by $\va_*(f)=\int(f)$.

  We check the axioms for the coproduct and counit.

  \begin{thm}\label{thm-hlh}
  $(H_{L^2}, *, e_K,\Delta_*,\varepsilon_*)$ is a Hopf algebra with
  approximate unit.\end{thm}
  \proof The proof consists of some lemmas. 
  For convenience we shall use the notation
  $$\Delta_*(f) = \sum_{(f)} f^{(1)} \otimes f^{(2)} = f^{1} \otimes
  f^{2},$$ where the sum on the right-hand side is an absolute convergent
  series in $H_{L^2} \hat{\otimes} H_{L^2}.$ \eee

  \begin{lem} \label{coact-axiom} The coproduct and the counit satisfy
  $$(\Delta_* \otimes id_{H_{L^2}}) \Delta_* = (id_{H_{L^2}} \otimes
  \Delta_*) \Delta_*,$$
  $$(\varepsilon_* \otimes id_{H_{L^2}})\Delta_* =
  (id_{H_{L^2}} \otimes \varepsilon_*)\Delta_* = id_{H_{L^2}},$$
 in the second equation we identify 
  ${\mathbb C} \hat{\otimes} H_{L^2}$ and $\HLh\hat{\otimes} {\mathbb C}$ with $ H_{L^2}$. \end{lem}
  \proof
  By the faithfulness of $\int$, the first equation is equivalent to
  $$\iiint (\Delta_* \otimes id)\Delta_*(f)\cdot  g\otimes h \otimes k =
  \iiint((id \otimes \Delta_*)\Delta_* f\cdot  g \otimes h\cdot k), \forall g,h,k \in
  H, f\in H_{L^2}$$
  By definition, the left hand side is equal to
  $$\iiint (f^{11} \otimes f^{12} \otimes f^2\cdot g\otimes h \otimes k) =
  \int(f^1\cdot gh)\int(f^2k) = \int fghk\cdot $$
  So is also the right hand side, whence we obtain the first equation of the
  lemma. 

  For the second equation of the lemma, we have, for all $g\in H$,
  $$\int(\varepsilon(f^1) f^2\cdot g = \int(f^1)\int(f^2g) = \int fg\cdot $$
  The lemma \ref{coact-axiom} is therefore proved. \eee

  \begin{lem} \label{copr-pr} The map
  $\Delta_* : H_{L^2} \to H_{L^2} \hat{\otimes} H_{L^2}$ is a homomorphism
  of algebras. \end{lem}
  \proof
  This is equivalent to the fact that for any $h\cdot k\in H$, $f,g\in H_{L^2}$,
  $$\int(f*g)hk = \int((f^1 * g^1)h)\cdot \int((f^2*g^2)k).$$
  The left hand side is equal to
  $$\int (fh_1k_1)\cdot \int(gh_2k_2) =
  \int(f^1h_1)\int(f^2k_1)\int(g^1h_2)\int(g^2k_2) =
  \int((f^1*g^1)h)\cdot \int((f^2*g^2)k).$$
  The Lemma \ref{copr-pr} is proved. \eee

  The above two lemmas imply that $H_{L^2}$ is a bialgebra with
  approximate unity. It remains to find an approximate antipode. By
  definition, we have to find a system $\{ S_{*K} \vert K \in C_0(\Lambda)
  \}$ of linear endomorphism on $H_{L^2}$, such that $m_*(S_{*K} \otimes \id)
  \Delta_*$ and $m_*(\id \otimes S_{*K})\Delta_*$ are approximate units.

  For any $h\in H_\lambda$, we define $S_*(h)\in H_{\lambda *}$ to be such that
  $\forall g \in H_\lambda$, 
 $$\int (hS(g)) = \int (S_*(h) g).$$ In fact, $S_*(h)$ can be
  computed explicitly
  $$S_*(h) = S(h_1)p(h_2)q(h_3),$$ where $q$ is the linear functional defined in
 Lemma \ref{lemssquare}, $p$ is given by
  $$p(h) = \int S(h_2) h_1, \forall h\in H.$$
  Set $S_{*\lambda}\vert_{H_\mu} = 0,\forall \mu \ne \lambda$.
  For any $K \subset C_0(\Lambda)$ set $S_{* K} = \bigoplus_{\lambda \in K}
  S_{* \lambda}$. 

  \begin{lem}
  $\{ S_{* K} \vert K\in C_0(\Lambda) \}$ is an approximate
  antipode.
  \end{lem}
  \proof
  We show that $$S_{* K}(f^1) * f^2 =\va_*(f)e_K= \int( f) e_K, \forall f \in H, K \in
  C_0(\Lambda).$$

 Since $\int S_{* K}(f) g = \int fS(g_K),$ we have
  $$\left.\begin{matrix}
  \int\left((S_{* K}(f^1) * f^2) g\right) & = & \int(S_{* K}(f^1)g_1) . \int (f^2 g_2)\\
     & = & \int(f^1S(g_{K1})).\int (f^2 g_2)\\
     & = & \int (f) \int (e_K g),
  \end{matrix}\right.$$
  whence the assertion follows. The lemma is proved. \eee

  We have therefore finished the proof of Theorem \ref{thm-hlh}. \eee

  \noindent{\bf Remark.}
  Since the antipode on $H$ is not involutive, the
  $*$-structure on $\check{H}$ defined by $h^\star = S(h^*)$ cannot be
  extended to $H_{L^2}$, because it is not continuous with respect to the
  given norm. 
  To overcome this obstruction we have to pass to the $C^*$ envelope of
  $H_{L^2}$ which is the subject of our next section.
 \section{The $C^*$-algebra $H_{C^*}$}\label{sec3}
 Let $H$ be a  compact Hopf $*$-algebra. Let $\star$ be the involutive map on \chH\ defined in the previous section $f^\star=S(f^*)$. 
 \begin{lem}\label{lemunitary-mod} All simple \chH-modules have the structure of $\star$-modules.\end{lem}
 \proof According to Theorem \ref{thmact-coact}, a simple \chH-comodule is equivalent to some module $\chV_\lambda$ induced from a simple unitary $H$-comodule $V_\lambda$. We check that, with respect to the given scalar product on $V_\lambda$,
 $$<h*v,w>=<v,h^\star*w>, \forall h\in \chH, v,w\in \chV_\lambda.$$
 Indeed, 
 \bbas <h*v,w>&=&<v_0,w>\overline{\int(hS(v_1))}\\
 &=&<v,w_0>\overline{\int(hw_1*)}\quad\mbox{ according to \rref{unitary-mod} }\\
 &=& <v,w_0>\int(w_1h^*)\\
 &=&<v,w_0>\int(S(h^*)S(w_1))\\
 &=&<v,h^\star*w>. \eqeee\eeas

  We introduce the following semi-norm on \chH\
 \begin{equation}\label{c*norm}
 \|f\|_{C^*}:=\sup_{\lambda\in\Lambda}\|\pi_\lambda(f)\|.\end{equation}\

 The lemma below will show that this is a norm, i.e., it is bounded.

 Let $\U=\U_\lambda$ be a unitary coefficient matrix corresponding to the simple unitary comodule $V_\lambda$ and $Q=Q_\lambda$ be the corresponding reflection matrix. If $0\neq f\in\chH_\lambda$, then $f=F_j^i\u_i^j$ for a complex matrix $F$. Since $f$ act on $\chV_\mu$ as zero if $\mu\neq \lambda$, then according to \rref{f-action},
 \begin{equation}\label{norm-f}
 \|f\|_{C^*}=\|\pi_\lambda(f)\|=\|FQ\|>0.\end{equation}

 \begin{lem}\label{leml2c*-norm}The $L^2$- and $C^*$-norms on \chH\ satisfy
 $$\|f\|_{L^2}\geq \|f\|_{C^*}.$$
 \end{lem}
 \proof 
 Since $\|f\|_{L^2}^2=\sum_\lambda\|f_\lambda\|_{L^2}^2$ and $\|f\|_{C^*}=\sup_\lambda\|f_\lambda\|_{C^*}$, it is sufficient to check the above equation for $f\in\chH_\lambda$. Thus $f=F^i_j\u_i^j$, for $(\u_i^j)=(\u_\lambda)_i^j$ a unitary multiplicative matrix.

  The desired inequality is then equivalent to
 $$\frac{\tr(FQF^*)}{\tr(Q)}\geq \frac{\|FQ\|^2}{\tr(Q)^2}.$$

 Since $Q$ is positive definite, $Q=T^*T$. Hence
 $$\|FQ\|^2=\|FTT^*\|\leq \|FT\|^2\|T^*\|^2\leq\tr(FQF^*)\tr(Q).$$
 Here we use the inequality
 $$\|A\|^2\leq \tr(AA^*). \eqeee$$

We define the algebra \Hcs\ to be the completion of \chH\ with respect to the norm $\|\cdot\|_{C^*}$. Since $\|\cdot\|_{C^*}$ is an operator norm, \Hcs\ is a $C^*$-algebra. By virtue of Lemma \ref{leml2c*-norm} we have 
 \bba\label{hl2hcs}\HLh\subset\Hcs.\eea

From the decomposition \rref{alg-decom} of \chH\ and the definition of the norm $\|\cdot\|$, we can easily obtain an isomorphism of the two $C^*$-algebras
\bba\label{gelfand}
\Hcs\cong {\prod_{\lambda\in\Lambda}}'\End_\bC(V_\lambda),\eea
where the product on the right-hand side of the equation is over all families $\{x_\lambda\in \End_\bC(V_\lambda)\}$ with $\|x_\lambda\|\lora 0$ as $\lambda\lora \infty$.
\begin{lem}\label{uni-repr} All irreducible unitary representations of \Hcs\ are finite-dimensional and irreducible over \chH.\end{lem}
\proof Let $\pi:\Hcs\lora \B(\H)$ be an irreducible unitary representation of \Hcs\ in a Hilbert space $\H$. Then $\H$ is a \chH-module. Let $M\subset\H$ be a simple \chH-submodule. Then it has a structure of a simple $H$-comodule, hence is finite dimensional over $\bC$. Thus $M$ is closed in $\H$.

Since \chH\ is dense in \Hcs\ and $M$ is closed in $\H$, $M$ is a representation of \Hcs. By the irreduciblity of $\pi$, we conclude that $M=\H$.\eee

As a corollary of Lemmas \ref{lemunitary-mod} and \ref{uni-repr}, we have
\begin{pro} There exists a 1-1 correspondence between irreducible unitary representations of \Hcs\ and those of \chH.\end{pro}

Thus \Hcs\ is a $C^*$-algebra of type I.

From the construction of \Hcs\, it is easy to see, that if $G$ is a compact group and $H={\BBb C}[G]$ -- the Hopf algebra of representative functions on $G$, which is a dense subalgebra of the algebra $C^\infty(G)$ of all continous complex valued functions on $G$, then $\Hcs\cong C^*(G)$ -- the group $C^*$-algebra of $G$.

In the case of a compact group $G$, the action of $G$ on any representation can be recovered by the action of $C^*(G)$, by using $\delta$-type sequences. In our case, we do not know, what the quantum group defining $H$ is. Therefore we have to introduce a Hopf algebra structure on \Hcs\ to exhibit the ``group'' property of this algebra.

 \begin{thm}\label{c*withau}	
 \Hcs\ is a Hopf $C^*$-algebra with approximate unit.\end{thm}
 \proof

 Recall that the topological tensor product $\Hcs\hat{\ot}\Hcs$ is the completion of $\chH\ot\chH$ with respect to the norm
 $$\|f\ot g\|_{C^*}:=\sup_{\lambda,\mu\in\Lambda}\|\pi_\lambda(f)\ot \pi_\mu(g)\|.$$

 For any $f\in\chH\subset\HLh$, $\Delta(f)\in\Hcs\hat{\ot}\Hcs.$ Thus, it is sufficient to show that
 $$\|\Delta(f)\|_{C^*}\leq\|f\|_{C^*}.$$
 But this inequality is obvious, as for any $\lambda,\mu\in\Lambda$, $V_\lambda\ot V_\mu$ is a representation of \HLh\ by means of the map $(\pi_\lambda\ot\pi_\mu)\Delta$, hence decomposes into a direct sum of irreducible representations, namely $\pi_\lambda\ot\pi_\mu\cong \oplus_\gamma c_{\lambda\mu}^\gamma \pi_\gamma$, consequently
 $$\|(\pi_\lambda\ot\pi_\mu)\Delta(f)\|\leq\sup_{\gamma\in\Lambda}\|\pi_\gamma(f)\|=\|f\|_{C^*}.$$

 To show that $\va_*$ extends onto \Hcs\ we remark that $\va_*:\chH\lora {\Bbb C}$ is multiplicative, i.e., induces a representation, for
 $$\va_*(f*g)=\int(f_1\int(f_2S(g)))=\int(f)\int(g)=\va_*(f)\va_*(g).$$
  Hence
 $$|\va_*(f)\leq \|f\|_{C^*}$$
 thus $\va_*$ extends onto \Hcs.

 According to Lemma \ref{lemunitary-mod}, for $f\in\chH$,
 $$\pi_\lambda(f^\star)=\pi_\lambda(f)^\star, \forall \lambda\in\Lambda.$$
 Therefore the involutive map $\star$ extends on $H_{C^*}$. 

 Using \rref{ssquare} and \rref{conv-pr2}, we can show that the unit element $e_\lambda$ in $H_\lambda$ satisfies $e_\lambda^\star=e_\lambda$. Since $\{e_K|K\in C_0(\Lambda)\}$ is an approximate unit of \chH\ and since \chH\ is dense in \Hcs, $\{e_K|K\in C_0(\Lambda)\}$ is an approximate unit on \Hcs. Consequently, $\{S_K|K\in C_0(\Lambda)\}$ is an approximate antipode on \Hcs.
 The proof of Theorem \ref{c*withau} is completed. \eee

\end{document}